\documentclass[12pt]{article}

\usepackage{amsmath}
\usepackage{amsfonts}

\newtheorem{theo}{Theorem}

\newtheorem{prop}{Proposition}
\newtheorem{lemm}{Lemma}
\newtheorem{defn}{Definition}

\newcommand{\lbl}{\label}

\newcommand{\eq}[1]{$(\ref{#1})$}

\def\N{\mathbb{N}}
\def\E{\mathbb{E}}
\def\0{{\bf 0}}
\def\Z{\mathbb{Z}}
\def\R{\mathbb{R}}

\def\Co{K}

\renewcommand{\E}{\mathbb E \,}

\newcommand{\testfn}{{\cal R}(Q_1^+)}

\newcommand{\LL}{{\cal L}}
\newcommand{\tB}{\tilde{B}}

\newcommand{\eqd}{\stackrel{{\cal D}}{=}}

\newcommand{\toas}{\stackrel{{{\rm a.s.}}}{\longrightarrow}}

\newcommand{\eqco}{\setcounter{equation}{0}}
\newcommand{\thco}{\setcounter{theo}{0}}
\newcommand{\prco}{\setcounter{prop}{0}}
\newcommand{\laco}{\setcounter{lemm}{0}}
\newcommand{\coco}{\setcounter{coro}{0}}
\newcommand{\cjco}{\setcounter{conj}{0}}

\newcommand{\deco}{\setcounter{defn}{0}}
\newcommand{\allco}{\eqco  \thco \prco \laco \coco \cjco \deco}
\newcommand{\Bx}{{\rm Box}}
\newcommand{\lad}{\la}
\newcommand{\lainf}{\la}
\newcommand{\qed}{\hfill{\rule[-.2mm]{3mm}{3mm}}}

\newcommand{\Po}{{\cal P}}
\newcommand{\X}{{\cal X}}
\newcommand{\A}{{\cal A}}
\newcommand{\Y}{{\cal Y}}

\def\la{{\lambda}}
\def\r{{\rho}}

\renewcommand{\P}{{{\cal P}}}

\newcommand{\Cov}{{\rm Cov}}

\newcommand{\Var}{{\rm Var}}

\newcommand{\diam}{{\rm diam}}

\newcommand{\XX}{{\cal X}}

\newcommand{\Cl}{{\rm Cl}}
\newcommand{\eps}{\varepsilon}
\def\bdm{\begin{displaymath}}
\newcommand{\edm}{\end{displaymath}}
\def\benu{\begin{enumerate}}
\def\eenu{\end{enumerate}}
\def\beqn{\begin{equation}}
\def\eeqn{\end{equation}}
\def\be{\begin{equation}}
\def\ee{\end{equation}}
\def\bea{\begin{eqnarray}}
\def\eea{\end{eqnarray}}
\newcommand{\bean}{\begin{eqnarray*}}
\newcommand{\eean}{\end{eqnarray*}}
\newcommand{\bear}{\begin{eqnarray}}
\newcommand{\eear}{\end{eqnarray}}

\renewcommand{\epsilon}{\varepsilon}

\def\R{\mathbb{R}}

\def\de{{\delta}}

\def\qed{\hfill\hbox{${\vcenter{\vbox{
    \hrule height 0.4pt\hbox{\vrule width 0.4pt height 6pt
    \kern5pt\vrule width 0.4pt}\hrule height 0.4pt}}}$}}

\def\la{{\lambda}}

\begin{document}

\title{\bf Gaussian limits for multidimensional 
random sequential packing at saturation (extended version)}

\author{T. Schreiber$^{1}$, Mathew D. Penrose and J. E. Yukich$^{2}$ }

\maketitle

\footnotetext{ {\em American Mathematical Society 2000 subject
classifications.} Primary 60F05, Secondary 60D05, 60K35} 

 \footnotetext{$~^1$
Research partially supported by the Polish Minister of Scientific
Research and Information Technology grant 1 P03A 018 28
(2005-2007)}

\footnotetext{$~^2$ Research supported in part by  NSF grant
DMS-0203720}

\begin{abstract}Consider the random sequential packing model with
infinite input and in any dimension.  When the input consists of
non-zero volume convex solids we
 show that the total number of solids accepted over
cubes of volume $\la$ is asymptotically normal as $\la \to
\infty$. We provide a rate of approximation to the normal and show
that the finite dimensional distributions  of the packing measures
converge to those of a mean zero generalized Gaussian field. The
method of proof involves showing that the collection of accepted
solids satisfies the weak spatial dependence condition known as
stabilization.
\end{abstract}

\section{Main results}\label{INTRO}

\allco

Given
 $d \in \N$ and
$\la \geq  1$,
 let $U_{1,\la},
U_{2,\la}, \ldots$  be a sequence of independent random $d$-vectors uniformly distributed
on the cube
$Q_\la := [0, \la^{1/d})^d$.
Let $S$ be a fixed bounded closed convex set in $\R^d$ with non-empty interior
(i.e., a  `solid')
 with
centroid at the origin $\0$ of $\R^d$ (for example, the unit ball),  and
for $i \in \N$,
let $S_{i,\la}$ be the  translate of $S$ with centroid
at $U_{i,\la}$. So ${\cal S}_{\la} :=
 (S_{i,\la})_{i \geq 1}$ is an infinite sequence of
 solids arriving at uniform random positions
in $Q_\la$  (the centroids lie in $Q_\la$ but the solids themselves
need not lie wholly inside $Q_\la$).

 Let the first solid   $S_{1,\la}$ be {\em
packed}, and recursively for $i=2,3, \ldots$, let
 the $i$-th solid
$S_{i,\la}$ be  packed if  it  does not overlap any solid
in $\{ S_{1,\la},\ldots,S_{i-1,\la} \}$ which has already been packed. If
not packed, the $i$-th solid is discarded;  we sometimes
use {\em accepted} as a synonym for `packed'.   This process, known as
{\em random sequential adsorption (RSA) with infinite input}, is
irreversible and terminates when it is not possible to accept
additional solids. At termination,  we say that the sequence of
solids  ${\cal S}_{\lainf}$
{\em jams} $Q_\la$ or {\em saturates} $Q_\la$. The {\em jamming} number
$N_{\la}:= N_{\la}({\cal S}_{\la})$ denotes the number
of solids accepted in $Q_\la$ at termination.
We use the words `jamming' and `saturation' interchangeably in
this paper.

Jamming numbers $N_{\la}$ arise naturally in the physical,
chemical, and biological sciences.  They are considered in the
description of the irreversible deposition of colloidal particles
on a substrate (see the survey \cite{ASZB} and the special volume
\cite{Pr}), hard core interactions (see the survey \cite{E}; also
\cite{To}),  adsorption modelling (see \cite{BP}
and the survey \cite{T})
and also in
the modelling of communication
and reservation protocols (see \cite{CFJ2, CFJP}).

The extensive body of experimental results related to the large
scale behavior of packing numbers stands in sharp contrast with
the limited collection of rigorous mathematical results,
especially in $d \geq 2$. The main obstacle to a rigorous
mathematical treatment of the packing process is that the short
range interactions of arriving particles create long range spatial
dependence, thus turning $N_{\la}$ into a sum of spatially
correlated random variables.

In the case where $d =1 $ and $S=[0,1]$,
  a famous result of R\'enyi \cite{Re} shows that {\em
jamming limit}, defined as  $\lim_{\la \to \infty} \la^{-1} \E
N_{\la} $, exists as an integral which evaluates to roughly
$0.748$;
also in this case, Mackenzie \cite{Mac} shows that
$\lim_{\la \to \infty} \la^{-1} \Var N_{\la}$
exists as an integral which evaluates to roughly $0.03815$.
Dvoretzky and Robbins \cite{DR} show
that the  jamming numbers $N_{\lad}$ are asymptotically
normal as $\la \to \infty$, but their  techniques do not
address the case $d >1$.

Since the above results were established in the 1960s,
progress in
extending them rigorously to higher dimensions
has been slow until recently.
Penrose \cite{Pe1}
 establishes
 the existence of a jamming limit for
any $d \geq 1$ and any choice of $S$, and also
 \cite{Pe2} obtains a CLT for a related model
(monolayer ballistic deposition with a rolling mechanism)
but comments in \cite{Pe2} that
 `Except in the case $d=1$ ...
a CLT for infinite-input continuum  RSA remains
elusive.'

In the present work  we show for any $d$ and $S$
that $\la^{-1} \Var N_{\la}$ converges to a
positive limit and that $N_\la$ satisfies a central limit theorem,
i.e.,
 the fluctuations of the random variable $N_{\la}$
are indeed Gaussian in the large $\la$ limit. This puts the recent
experimental results and Monte Carlo simulations of Quintanilla
and Torquato \cite{QT} and Torquato (ch. 11.4 of \cite{To}) on
rigorous footing. We also provide a bound on the rate of
convergence to the normal, and on the rate of convergence of
$\la^{-1} \E N_\la$ to the jamming limit.

Throughout ${\cal N}(0, 1)$ denotes a mean zero normal random
variable with variance one.

\begin{theo} \label{CLT1}
Let ${\cal S}_{\lainf}$ be as above and put $N_{\la} :=
N_{\lad}({\cal S}_{\la})$. There are constants  $\mu:=\mu(S,d) \in
(0, \infty)$ and $\sigma^2 := \sigma^2(S,d) \in (0, \infty)$ such
 that as $\la \to \infty$ we have  \bea
| \la^{-1} \E N_{\la} -  \mu | = O(\la^{-1/d})
\lbl{meanrate}
\eea
and
$
\la^{-1} \Var N_{\la} \to \sigma^2
$
with
\bea
\sup_{t \in \R}   \left|  P \left[  { {N_{\lad} - \E N_{\lad}}   \over (
\Var N_{\lad} )^{1/2} }  \leq t  \right]
 - P[ {\cal N}(0, 1) \leq t] \right|
  = O( (\log \la)^{3d} \la^{-1/2}).
\lbl{jamrate}
\eea
\end{theo}
\vskip.5cm

The process of accepted solids in $Q_\la$ induces a natural random
point measure  $\nu_\la$ on $[0,1]^d$ given by \be
 \label{defnu} \nu_\la := \sum_{i = 1}^{\infty}
\delta_{\la^{-1/d} U_{i, \la}} {\bf 1}_{\{S_{i,\la} \mbox{is accepted} \}}
\ee
where $\de_x$ stands for the unit point
mass at $x$.
It also induces a natural random volume measure $\nu'_\la$ on $\R^d$,
normalized to have the same total measure as $\nu_\la$,
defined for all Borel $A \subseteq \R^d$ by
\be
 \label{defnu2} \nu'_\la (A) := \frac{\la}{|S|} \left|
A \cap \left( \bigcup[
 \la^{-1/d} S_{i,\la} :
i \geq 1, S_{i,\la} \mbox{ is accepted} ] \right) \right| \ee
where $|\cdot|$
denotes Lebesgue measure and
 $\la^{-1/d}A :=\{\la^{-1/d}x: x \in A\}$.
The measure $\nu'_\la$ is not necessarily
supported by $Q_1$ due to boundary effects,
but for $\la >1$ it is supported by $Q_1^+$, where we set
$Q_1^+ := [-1,2)^d$ (a fattened version of $Q_1$).

 Let $\bar\nu_\la := \nu_\la - {\Bbb E}[\nu_\la]$
 and $\bar\nu'_\la := \nu'_\la - {\Bbb E}[\nu'_\la]$.
Let
$\testfn$
 denote the class of
 bounded,  almost everywhere continuous functions on
$Q_1^+$.
 For $f \in \testfn $ and $\mu$ a
 signed measure on $\R^d$ with finite total mass, let
$\langle f, \mu \rangle := \int_{\R} f d\mu$.
The following theorem provides the limit theory (law of large numbers
and central limit theorems)  for the integrals
of test functions $f \in \testfn$ against the
 random point measure  $\nu_\la$ and the random volume
measure $\nu'_\la$ induced by the packing process.
In particular, it
%
shows that the  finite dimensional
distributions of the centered packing point measures $(\bar\nu_\la)_\la
$ converge to those of a
certain mean zero generalized Gaussian field,
namely white noise on $Q_1$ with variance $\sigma^2$
per unit volume,
and likewise for the centered packing volume measures
$(\bar\nu'_\la)_\la$.


\begin{theo} \label{CLT2}
Let $\mu$ and $\sigma^2 $ be as in Theorem \ref{CLT1}. Then
for any $f, g$ in $\testfn$,
$$
\lim_{\la \to \infty} \la^{-1}
 \E
[ \langle f, \nu_\la \rangle] =  \mu \int_{[0,1]^d} f(x) dx
$$
and
$$
\lim_{\la \to \infty} \la^{-1}
 \Cov(
\langle f, \nu_\la \rangle,
\langle g, \nu_\la \rangle )
=
 \sigma^2 \int_{[0,1]^d} f(x) g(x) dx.
$$
Also,
 the finite-dimensional distributions of
the random field
  $(\la^{-1/2} \langle  f, \bar \nu_\la \rangle ,f \in \testfn)$
 converge
as $\la \to \infty$  to those of a mean zero generalized Gaussian
field with covariance kernel
$$
 (f,g) \mapsto
 \sigma^2 \int_{[0,1]^d} f(x) g(x) dx,\;\;
f,g \in  \testfn.
  $$
Moreover, the same conclusions hold with $\nu_\la$
and $\bar{\nu}_\la$ replaced by $\nu'_\la$ and $\bar{\nu}'_\la$
respectively.
\end{theo}

\vskip.5cm



{\bf Remarks.}

1.  {\em Finite input}.
Let $\tau \in (0,\infty)$ and let
 $\lceil x \rceil$ denote the smallest integer greater than
 or equal to $x$.
Inputting only the first $\lceil \la \tau \rceil$
 solids of the sequence ${\cal S}_{\lainf}$
 yields RSA packing of the cube $Q_\la$ with {\em
finite input}.
The finite-input
packing number, i.e., the total
number of solids accepted from $S_{1,\la},S_{2,\la},..., S_{\lceil
\tau \la \rceil, \la}$,
  is asymptotically normal as $\la \to \infty$ with $\tau $ fixed.  This is
  proved in \cite{PY2},
and extended in
\cite{BY2} to the case
 where the spatial coordinates come from a
non-homogeneous point process. Packing measures induced by RSA
packing with finite input have finite dimensional distributions
converging to those of a mean zero generalized Gaussian field with
a covariance structure depending upon the underlying density of
points \cite{BY2}.

 2. {\em Stabilization.}
One might expect that the restriction of the packing measure
$\nu_\la $ or $\nu'_\la$ to a localized region of space
depends only on
incoming particles
with `nearby' spatial locations, in some
well-defined sense.
This local dependency property is denoted {\em stabilization};
when the region of spatial
dependency has a diameter with an exponentially decaying tail, it
is called {\em exponential stabilization}. These notions
are spelt out in general terms in Section 2.
 Theorem \ref{basethm1}
provides a general spatial limit theory for
exponentially stabilizing measures;
 this is  an infinite-input
analog to known  results
\cite{BY2,Pe4,Pe5,Pe6} for
the finite-input setting, and is of independent interest.

A form of stabilization for infinite input RSA  was proved in
\cite{Pe1}, but without any tail bounds.  {\em Exponential}
stabilization in the infinite input setting is perhaps not
surprising, but it has been challenging to rigorously establish
this key localization feature.  In Section \ref{secstab}, we show
that  infinite-input packing measures
 stabilize  exponentially, so that the general results
of Section \ref{secstein} are applicable to these measures.

3. {\em Related  models} in the literature (see e.g. \cite{PY2})
include cooperative sequential adsorption, RSA with solids of
random size or shape, ballistic deposition with a rolling
mechanism, and spatial birth-growth models.
For all of these models, limit theorems in the finite-input setting
are discussed in \cite{Pe2}. It seems likely
that these can be
extended to the infinite-input setting using the methods
of this paper, although we  do not discuss any of them in detail.
Nor do we consider  non-homogeneous
point processes as input.


4. {\em Rates of convergence}.  Even in $d = 1$, the rate given by
Theorem \ref{CLT1} is new. Quintanilla and Torquato \cite{QT} use
Monte Carlo simulations to predict convergence of the distribution
function for $N_{\lad}$ to that of a normal, but they do not
obtain rates.  Penrose and Yukich \cite{PY5} obtain rates of
approximation to the normal for RSA packing with finite (Poisson)
 input.

5. {\em Numerical values.} We do not provide any new analytical
methods for
 computing numerical
values of $\mu$ and $\sigma^2$ when $d \geq 1$.


6. {\em Jamming variability.}
A significant amount of work is needed (see
Section \ref{secJV})
to show that the
limiting variance $\sigma^2$ in Theorems \ref{CLT1} and
\ref{CLT2} is non-zero,
and we prove this using the following  notions.

 Given
 $L > 0$, we shall us say that a   point set
$\eta \subset \R^d \setminus [0,L]^d$ is {\em admissible} if the
translates of $S$ centered at the points of $\eta$ are
non-overlapping. Given such an $\eta$, let $N[[0,L]^d | \eta]$
denote the (random) number of solids from the sequence
${\cal S}_{L^d}$ which are packed in $[0,L]^d$
 given the
 {\em pre-packed configuration $\eta$}.
 In other words,
$N[[0,L]^d | \eta]$ arises as the number of solids packed in
$[0,L]^d$ in the course of the usual infinite input packing
process subject to the additional rule that an incoming solid is
discarded should it overlap any  solid centered at a point of
$\eta$.
  Say that the convex body
$S$ has {\em jamming variability} if there exists a $L > 0$ such
that $\inf_{\eta} \Var N[[0,L]^d | \eta] > 0$
with the infimum taken over admissible point sets $\eta \subset
\R^d \setminus [0,L]^d$.

In Proposition \ref{CLT3} we shall show that
 each  bounded convex body $S \subset \R^d$ with non-empty interior
 has jamming variability.

7.
We let $d_S$
stand for the diameter of $S$.
In our proofs,
{\em we shall assume that $2 d_S < 1$}.
This assumption entails
no loss  of generality, since once we have proved Theorems
\ref{CLT1} and \ref{CLT2}
under this assumption, the results follow for general $S$ by
obvious scaling arguments.

\section{Terminology, auxiliary results}\label{TEAUX}
\label{secstein}
\allco
Let $\R_+ := [0,\infty)$.
 Given a point $(x_1,\ldots,x_d,t) = (x,t) \in \R^d  \times \R_+$,
the first $d$ coordinates of the point
 will be interpreted as spatial
  components with the $(d+1)$-st regarded as a time mark.
  Let us say a point set $\X \subset \R^d \times \R_+$ is
  {\em temporally locally finite} (or TLF for short)
 if $\X \cap (\R^d \times [0,t])$
  is finite for all $t >0$. Loosely speaking, $\X$ is TLF
if it is finite in the spatial  directions and
locally finite in the time direction.

In this section we  adapt the general results and
terminology from \cite{BY2,Pe5,Pe6,PY5} on
limit theory for stabilizing
spatial measures defined in terms of finite point sets in $\R^d$, to
to the setting of spatial measures
defined in terms of  TLF point sets in $\R^d \times \R_+$ (typically
obtained as Poisson processes).
In subsequent sections, we show that these general results
can be applied to obtain  the limit theorems for RSA
described in Section \ref{INTRO}.

  For $x \in \R^d$ and $r > 0$, let
$B_r(x)$ denote the Euclidean ball centered at $x$ of radius $r$.
We abbreviate $B_r(\0)$ by $B_r$.
Given $\X \subset \R^d \times \R_+$,  $a > 0$ and $y \in \R^d$, we
let $y+ a\X:= \{(y+ax,t): (x,t) \in \X\}$;
 in other words, scalar multiplication and
translation  on $\R^d \times \R_+$ \  act only on the spatial
components. For $A \subset \R^d$ we write $y+aA $ for $\{y+ax: x \in A\}$;
also, we
write $\partial A$ for the boundary of $A$, and
  write $A_+$ for $A \times \R_+$.
For nonempty subsets $A,A'$ of $\R^d$,  write $D_2(A,A')$
for the Euclidean distance between them, i.e. $D_2(A,A')
:= \inf\{|x-y|:x \in A, y \in A'\}$.

Let $\xi(\X,A)$ be an
 $\R_+$-valued function  defined for all pairs $(\X,A)$,
where $\X$ is a TLF subset of $\R^d \times \R_+$ and $A$ is a
Borel subset of $\R^d$. Throughout this section we make the
following assumptions on $\xi$:
\begin{enumerate}
\item
  $\xi(\cdot,A)$ is   measurable
for each Borel $A$,
\item
 $\xi(\X,\cdot)$ is a finite measure
on $\R^d$ for each TLF  $\XX \subset \R^d \times \R_+$,
\item
 $\xi$ is {\em translation invariant}, that is
$\xi(i + \X, i +A) = \xi (\X, A)$ for all $i \in \Z^d$, all
TLF
$\XX \subset \R^d \times \R_+$, and all Borel $A \subseteq \R^d$,
\item
 $\xi$ is {\em uniformly locally bounded} (or just {\em bounded}
for short) in the sense
 that there is a finite constant
$||\xi||_{\infty}$ such that
for all TLF $\X \subset \R^d  \times \R^+$ we have
\bea
\xi (\X,[0,1]^d) \leq ||\xi||_{\infty}.
\lbl{bounded}
\eea
\item $\xi$ is {\em locally supported}, i.e. there exists a
 constant $\rho $ such that $\xi(\X,A) =0$ whenever $D_2(\X,A)
> \rho$.
\end{enumerate}
Note that if $\xi(\X,\cdot)$ is a point measure supported by the
points of $\X$, then $\xi$ is locally supported (in fact, in this
case we can set $\rho=0$).

For all $\la > 0$, let $\P_{\la}$ denote a homogeneous Poisson
point process in $\R^d \times \R_+$ with intensity measure $\la dx
\times ds$, with $dx$ denoting Lebesgue measure on $\R^d$ and $ds$
Lebesgue measure on $\R_+$.  We put $\P := \P_1$.

%
%
%
Thermodynamic limits and central limit theorems for functionals in
geometric probability are often proved by showing that the
functionals satisfy a type of local spatial dependence known as
stabilization
\cite{BY2, Pe4, Pe5, Pe6, PY2, PY4, SY1}
and that
will be our goal here as well. First, we adapt the definitions in
\cite{BY2, Pe4, Pe5}
 to the context of measures defined in
terms of  TLF point sets
in $\R^d$. Recall that $Q_\la $ denotes the cube $[0,\la^{1/d})^d$.

\begin{defn}\label{stab}  We say $\xi$ is {\em homogeneously
 stabilizing}
if there exists an a.s. finite  random variable $R'$ 
(a {\em radius of homogeneous stabilization} for $\xi$) such that
for all TLF $\X \subset (\R^d \setminus B_{R'})_+$
we have
\be
\label{stab0}
\xi ((\Po \cap (B_{R'})_+) \cup \X,Q_1) =
\xi (\Po \cap (B_{R'})_+,Q_1).
\ee
We say $\xi $ is
 {\em exponentially stabilizing}
 if (i) it is homogeneously stabilizing and $R'$ can be chosen so
that
$ \limsup_{ L \to \infty} L^{-1} \log P[R'>L] < 0$, and
(ii) for all $\la \geq 1 $ and all
  $i \in   \Z^d$,
there exists a random
variable $R:=R^{\xi}(i,\la)$ (a {\em radius of stabilization} for
$\xi$ at $i$ with respect to $\Po $ in $(Q_\la)_+)$ such that for
all TLF
 ${\X } \subset \left[ Q_\la \ \setminus B_{ R} (i)
\right]_+$, and all Borel $A \subseteq Q_1$,
 we have
 \be
  \label{stab1}
 \xi \left(  (\P \cap [B_{ R}(i) \cap Q_\la]_+) \cup \X ,i+ A
 \right)
=
  \xi \left(    \P  \cap [B_{ R}(i) \cap Q_\la]_+ ,i+ A
  \right)
\ee
and moreover
the tail probability $\tau(L)$ defined for $L
>0$ by
\bea \tau(L) := \sup_{\la \geq 1,\ i \in \Z^d
}
P[R^{\xi}(i,\la)
> L] \label{taudef} \eea satisfies \ \
$ \limsup_{ L \to \infty} L^{-1} \log \tau(L) < 0.
 \label{tau1}$
\end{defn}

Loosely speaking, $R:=R^{\xi}(i,\la)$ is a radius of stabilization
if the  $\xi$-measure on $i+Q_1$
is unaffected by changes to the Poisson  points
outside $B_{R}(i)$ (but inside $Q_\la$).
When $\xi$ is homogeneously stabilizing, the limit
$$
\xi( \P,i+ Q_1 ):= \lim_{r \to \infty} \xi \left(\P
\cap (B_{ r}(i))_+ ,i+Q_1
\right)
$$
exists almost surely for all $i \in \Z^d$. The random variables
$(\xi(\P, i + Q_1), i \in \Z^d)$ form a stationary random field.

Given
$\xi$,
 for all $\la > 0$, all TLF $\X \subset \R^d \times \R_+$, and
all Borel $A \subset \R^d$
we let
$
\xi_{\la}(\X, A) := \xi (\la^{1/d}\X,\la^{1/d} A  ).
$
 Define the random measure $\mu_\la^\xi$ on
 $\R^d $  by
 \be
\label{defmu} {\mu}_{\la }^{\xi}( \ \cdot \ ):=
 \xi_\la( {\cal P}_{\la} \cap Q_1 ,  \cdot)
\ee
and the centered  version $\overline{\mu}^{\xi}_{\la } :=
{\mu}_{\la }^{\xi} - \E[{\mu}_{\la }^{\xi}]$.
By the assumed locally supported property of $\xi$,
$\mu_\la$ is supported by the fattened cube $Q_1^+ := [-1,2)^d$
for large enough $\la$.

If $\xi$ is stabilizing,
define
$\mu(\xi) := \E [\xi(\P,Q_1 )] $ and
and if $\xi$ is exponentially stabilizing, define
$$
 \sigma^2(\xi) :=
    \sum_{i \in \Z^d } \Cov \left[ \xi(\P,Q_1),
 \ \xi(\P,i+Q_1) \right],
 $$
where the sum can be shown to converge absolutely
by exponential stabilization and \eq{bounded}.
The following general theorem provides  laws of large numbers and
 normal approximation results for $\langle f, \mu_\la^\xi\rangle$,
 suitably scaled and centered, for $f \in \testfn$.
This  set of results for measures determined by TLF point sets
is similar to previously known results
  for measures determined by finite point sets
 (Theorem 2.1 of \cite{PY4},
Theorem 2.1 of \cite{BY2},
 Theorem 2.3 of \cite{BY2},
and
 Corollary 2.4 of \cite{PY5}).


\begin{theo}\label{basethm1}
 Suppose that
 $\xi$ is
  exponentially stabilizing. Then as $\la \to \infty$,
  for $f$ and $g$ in $\testfn$ we have
\be
 \label{explim2}
 \lim_{\la \to \infty} \la^{-1} \E [\langle f,
   \mu_\la^\xi \rangle]
=  \mu(\xi) \int_{[0,1]^d}
   f(x)dx
   \ee
   and
    \be
  \label{varlim2} \lim_{\la \to \infty} \la^{-1} \Cov
  [\langle f, \mu_\la^\xi \rangle,
 \langle g, \mu_\la^\xi \rangle]
  = \sigma^2(\xi) \int_{[0,1]^d}
   f(x)g(x)dx. \ee
Also,
\be
 \label{explim}
| \la^{-1} \E
 [\mu_{\la}^{\xi}(Q_1^+)]  -  \mu(\xi) | = O( \la^{-1/d} ).
 \ee
 Moreover, if $\sigma^2(\xi) > 0$ then
\be \label{rate}
 \sup_{t \in \R}  \left| P \left[ { \mu^{\xi}_\la(Q_1^+) - \E
 [\mu^{\xi}_\la(Q_1^+)] \over ( \Var [\mu^{\xi}_\la(Q_1^+)] )^{1/2}  }
 \leq t \right] - P[{\cal N}(0,1) \leq t]
\right|
= O(  (\log \la )^{3d} \la^{-1/2})
\ee
%
and the finite-dimensional distributions of the random field
    $(\la^{-1/2} \langle  f, \bar \mu_{\la}^{\xi} \rangle ,f \in \testfn)$ converge
as $\la \to \infty$  to those of a mean zero generalized Gaussian
field with covariance kernel
%
 $$
  (f,g) \mapsto \sigma^2(\xi) \int_{[0,1]^d} f(x) g(x) dx,\; \ \ f,g \in
\testfn.
 $$
\end{theo}

\vskip.5cm

We shall use Theorem \ref{basethm1}
 to prove the results
on RSA described in Section \ref{INTRO}.
It seems likely that Theorem \ref{basethm1} can
also be applied
to obtain similar results  for the related  models
listed in Remark 3 of Section \ref{INTRO}.
For some of these, certain generalizations of
Theorem \ref{basethm1} may be needed; for example,
in some cases one may need
to allow for the Poisson points to carry independent identically
distributed random marks, and in others the boundedness condition
 \eq{bounded} may need to be
relaxed to a moments condition.
It seems likely that little change to the   proof of
 Theorem \ref{basethm1}  will be needed to cover
 these generalizations.

As we shall see shortly, the thermodynamic limits (\ref{explim2})
and (\ref{explim})  do not require exponential decay of the
stabilization radius for $\xi$, but in fact hold under weaker
decay conditions. We expect that
 (\ref{varlim2}) also holds under weaker decay conditions
on the stabilization radius,
and also that the boundedness condition \eq{bounded} can be
relaxed to a moments condition in Theorem \ref{basethm1},
but for simplicity we shall assume throughout
that $\xi$ is exponentially stabilizing and satisfies \eq{bounded}.
Also, if we restrict attention to $f$ supported by $Q_1$,
we do not need the condition that $\xi$ be locally supported.

\vskip.5cm

The rest of this section is devoted to proving Theorem \ref{basethm1}.
We shall use the following notation.
 Given $f \in \testfn$,
we extend $f$ to the whole of $\R^d$ by setting
 $f(x)=0$ for  $x \in \R^d \setminus Q_1^+$.
Given TLF $\X \subset (\R^d)_+$, and $\la >0$, write $\langle f,
\xi_\la(\X)\rangle$ for $\int_{\R^d} f(x) \xi_\la(\X,dx)$
 (the integral of $f$
with respect to the measure $\xi_\la(\X,\cdot)$). For
  $j \in \la^{-1/d} \Z^d $, let
 $f_{\la,j}: \R^d \to \R$  be given by $f_{\la,j}(x) =f(x)$
for $x \in j + Q_{1/\la}$, and $f_{\la,j}(x) =0$ otherwise. Then
 \bea \langle f,\mu_\la^{\xi}
\rangle = \sum_{j \in \la^{-1/d} \Z^d} \langle  f_{\la,
j},\mu_\la^{\xi} \rangle . \lbl{0123} \eea Also, let \bean
\overline{f}(\la,j)
 := \sup \{ f(x): x \in j + Q_{1/\la } \} ; ~~~~~
\underline{f}(\la,j) := \inf \{ f(x): x \in j + Q_{1/\la } \}.
\eean
For $x \in \R^d$ let $i_\la(x)$ be the choice of
$i \in \la^{-1/d} \Z^d$ such that
$x \in i + Q_{1/\la} $. \\

{\em Proof of  (\ref{explim2}).}
Let $f \in \testfn$. Then by \eq{0123}, we have \bea \la^{-1} \E
[\langle f, \mu_\la^{\xi} \rangle] = \la^{-1} \sum_{j \in
\la^{-1/d} \Z^d } \E [\langle f_{\la,j}, \xi_\la(\P_\la \cap
(Q_1)_+) \rangle]
\nonumber \\
= \int_{\R^d} \E[
 \langle f_{\la,i_\la(x)}, \xi_\la(\P_\la \cap (Q_1)_+) \rangle ] d x.
\lbl{0116d}
\eea

For $x \in \R^d \setminus \partial Q_1$,
 with $f$ continuous at $x$, we assert that
as $\la \to \infty$, \bea \E [ \langle f_{\la,i_\la(x)},
\xi_\la(\P_\la \cap (Q_1)_+) \rangle] \to \mu(\xi) f(x)1_{Q_1}(x).
\lbl{0116a} \eea This clearly holds for $x \in \R^d \setminus
[0,1]^d$, since both sides are zero for large $\la$, by the
locally supported property of $\xi$. To see \eq{0116a} for $x \in
(0,1)^d$, observe that
the left side
 has the upper bound
\bea \E  [\langle f_{\la,i_\la(x)}, \xi_\la(\P_\la \cap (Q_1)_+)
\rangle] \leq \overline{ f}(\la,i_\la(x))
 \E [ \xi_\la(  \Po_\la \cap (Q_1)_+ ,  i_\la( x) +Q_{1/\la} ) ]
\nonumber \\
=
\overline{ f}(\la,i_\la(x))
 \E [ \xi(  \Po \cap (Q_\la)_+ ,  i_1(\la^{1/d} x) +Q_1 ) ],
\lbl{0116b}
\eea
and has a similar
lower bound with
$\underline{ f}(\la,i_\la(x))$
instead of
$\overline{ f}(\la,i_\la(x))$.
If $f$ is continuous at $x$, then both
$\underline{ f}(\la,i_\la(x))$ and $\overline{ f}(\la,i_\la(x))$
tend to $f(x)$,  so to prove \eq{0116a} it suffices to show the
expectation in the last line of \eq{0116b} converges to
$\mu(\xi)$. By translation invariance, this expectation equals
\bean
 \E [
 \xi(  \Po \cap ( -i_1(\la^{1/d} x) + Q_\la)_+ , Q _1 )].
\eean
For $x$  in the interior of $Q_1$,
the set
 $ -i_1(\la^{1/d} x) + Q_\la$ has limit set
 $\R^d$ as $\la \to \infty$, i.e. for any $r < \infty$ the  ball $B_r$
is contained in
 $ -i_1(\la^{1/d} x) + Q_\la$ for large enough $\la$.
Hence by stabilization,
\bea
 \xi(  \Po \cap ( -i_1(\la^{1/d} x) + Q_\la)_+ , Q_1 )
\toas \xi( \Po, Q_1 )
\lbl{0116c}
\eea
and by \eq{bounded},
 the corresponding expectations converge.
This demonstrates \eq{0116a}.

The integrand in \eq{0116d} is dominated by  a constant for
$x \in Q_1^+$,
and is zero  for $x \notin Q_1^+$.
So by \eq{0116a}
and dominated convergence applied to \eq{0116d}, we obtain \eq{explim2}.
$\qed$ \\

{\em Proof of \eq{explim}.}
For this proof,
 set $f(x) \equiv 1$ on $Q_1^+$.
 We need to bound the error term in  \eq{explim2} for this choice of
$f$, which we do by using
\eq{0116d} again.
For $x \in \R^d$, let $X(x,\la)$ be the integrand in \eq{0116d}, i.e.  set
$X(x,\la) :=  \langle f_{\la,i_\la(x)}, \xi_\la(\P_\la \cap (Q_1)_+) \rangle$
with our current choice of $f$; also set
 $Y(x,\la) :=  \xi ( \Po \cap (-i_1(\la^{1/d} x) + Q_\la)_+,  Q_1 )$.
If  $x \in (0,1-\la^{-1/d})^d$, then
\bea
 \E[X(x,\la)]
= \E [ \xi_\la ( \Po_\la \cap (Q_1)_+ , i_\la(x) + Q_{1/\la} ) ]
\nonumber \\
= \E [ \xi ( \Po \cap (Q_\la)_+, i_1(\la^{1/d} x) + Q_1 )]
\nonumber \\
= \E [ \xi ( \Po \cap (-i_1(\la^{1/d} x) + Q_\la)_+,  Q_1 )] = \E[Y(x,\la)].
~~~
\lbl{0116e} \eea
Abbreviating
 the Euclidean
distance
 $D_2(\{y\},A)$
by
 $D_2(y,A)$ is
we have
\bean
D_2( \0, \partial( -i_1(\la^{1/d} x) + Q_\la) )
=
D_2( i_1(\la^{1/d} x)  ,\partial Q_\la )
\\
\geq
D_2(\la^{1/d} x,\partial Q_\la)
  -\sqrt{d}
=
 \la^{1/d} D_2(x, \partial Q_1)
-\sqrt{d} .
\eean
Hence the ball
$B_{
\la^{1/d} D_2(x,\partial Q_1) - \sqrt{d}}$
is contained in the box  $-i_1(\la^{1/d} x) + Q_\la$,
so with $R'$ denoting the radius of homogeneous stabilization
of $\xi$,
\bea
 Y(x,\la) {\bf 1}\{ R' <  \la^{1/d}D_2(x,\partial Q_1) - \sqrt{d}\}
= \xi ( \Po , Q_1)
  {\bf 1}\{ R' <  \la^{1/d} D_2(x,\partial Q_1) - \sqrt{d}\} .
\nonumber \\
\lbl{0117g} \eea Set $\mu := \mu(\xi) = \E [\xi (\Po,Q_1)]$. By
\eq{0116e}, \eq{0117g} and \eq{bounded}, we have
  for $x \in (0,1-\la^{-1/d})^d$ that
\bean
|\E[X(x,\la)] - \mu| =
|\E[Y(x,\la)] - \mu|
\\
=
|\E[(Y(x,\la) - \xi(\Po, Q_1) )
  {\bf 1}\{ R' \geq
 \la^{1/d} D_2(x,\partial Q_1) - \sqrt{d}\} ]|
\\
  \leq 2 \|\xi\|_\infty
P[ R' >
 \la^{1/d} D_2(x, \partial Q_1) - \sqrt{d} ] \eean
and so by exponential stabilization,
 there is a constant $K>0$ such that
\bea
|\E[X(x,\la)] - \mu|
\leq K \exp( -  \la^{1/d} D_2(x,\partial Q_1) /K).
\lbl{0116f}
\eea
Also by \eq{bounded}, for suitable $K$ the same bound
\eq{0116f}  for
 holds trivially for
\linebreak
$x \in Q_1 \setminus  (0,1-\la^{-1/d})^d$, and hence
\eq{0116f} holds for
all $x \in Q_1$.
By \eq{0116f}, it is straightforward to deduce that
\bea
\int_{Q_1} | \E[X(x,\la)] - \mu  | dx = O(\la^{-1/d} ).
\lbl{0123a}
\eea
Also, for $x \in \R^d \setminus Q_1$ with $D_2(x,\partial Q_1) > \la^{-1/d}$
we have $\E[X(x,\la)]=0$, and $X(x,\la)$ is uniformly bounded by
\eq{bounded}, so that
$$
\int_{\R^d \setminus Q_1} | \E[X(x,\la)]   | dx = O(\la^{-1/d} ).
$$
Combining this with \eq{0123a} and using \eq{0116d} gives us
\eq{explim}. $\qed$ \\


{\em Proof of \eq{varlim2}.}
 Let $f \in \testfn$ and assume $f$ is nonnegative.
By linearity, it
 suffices to prove \eq{varlim2} in the case
where $f$ is nonnegative  and $f\equiv g$, so we now assume this.
First, we assert that there is a constant $K$, independent of $\la$, such that
for all $\la \geq 1$ and all $i \in \la^{-1/d} \Z^d$,
$z \in \Z^d$, we have
\bea |
\Cov [ \langle
f_{\la,i},
\xi_\la(\Po \cap (Q_1)_+ ) \rangle,
\langle
 f_{\la,i+ \la^{-1/d}z} ,
\xi_\la(\Po \cap (Q_1)_+)
 \rangle ] |
\leq K \exp( - |z| /K). \lbl{expcov} \eea This can be proved by
arguments similar to those in,
 e.g.,  the proof  of Lemma 4.1 in \cite{BY2} or that of Lemma 4.2
in \cite{Pe6}.
By \eq{0123}, we have
\bea
\la^{-1}  \Var [\langle f,
   \mu^{\xi}_\la \rangle]  =
 \la^{-1}
 \sum_{i , j \in \la^{-1/d} \Z^d }
 \Cov[
\langle
f_{\la,i},
\xi_{\la}(\P_\la \cap (Q_1)_+) \rangle,
        \langle f_{\la,j}, \xi_\la( \P_\la \cap (Q_1)_+) \rangle]
%
%
%
%
\nonumber \\
= \int_{\R^d} dx \sum_{z \in \Z^d } \Cov [ \langle
 f_{\la,i_\la(x)},
\xi(\Po_\la \cap (Q_1)_+)
\rangle,
\langle
 f_{\la,i_\la(x)+\la^{-1/d}z},
\xi(\Po_\la \cap (Q_1)_+)
 \rangle )]
\nonumber \\
\lbl{0116g} \eea where the inner sum converges absolutely by
\eq{expcov} and is zero for $x \notin Q_1^+$.

Fix $x  \in (0,1)^d$ and $z \in \Z^d$, with $f$ continuous at $x$.
Then  we have the upper bound \bea \E [\langle
 f_{\la,i_\la(x)} ,
 \xi_\la(\Po_\la \cap (Q_1)_+)
\rangle
\langle
 f_{\la,i_\la(x)+\la^{-1/d}z} ,
\xi_\la(\Po_\la \cap (Q_1)_+)
 \rangle]
\nonumber \\
\leq
\overline{f}(i_\la(x),\la) \,
\overline{f}(i_\la(x)+ \la^{-1/d}z,\la )
\nonumber \\
 \times \E
[\xi_\la(\Po_\la \cap (Q_1)_+, i_\la(x) + Q_{1/\la}) \xi_\la(
\Po_\la \cap (Q_1)_+, i_\la(x)+ \la^{-1/d}z + Q_{1/\la})]
 \lbl{0117c} \eea 
and a similar  lower bound with
$\overline{f}(i_\la(x),\la) \: \overline{f}(i_\la(x)+ \la^{-1/d}z,\la) $
replaced by
\linebreak
$\underline{f}(i_\la(x),\la) \: \underline{f}(i_\la(x)+ \la^{-1/d}z,\la)  $.
  Note that both
$\overline{f}(i_\la(x),\la) \; \overline{f}(i_\la(x)+
\la^{-1/d}z,\la) $ and $\underline{f}(i_\la(x),\la)\:
\underline{f}(i_\la(x)+ \la^{-1/d}z,\la) $ converge as $\la \to
\infty$ to $f^2(x)$.

By scaling and translation invariance of $\xi$, we have  \bean \E
[ \xi_\la(\Po_\la \cap (Q_1)_+, i_\la(x) + Q_{1/\la}) \xi_\la(
\Po_\la \cap (Q_1)_+, i_\la(x)+ \la^{-1/d}z + Q_{1/\la}) ]
\\
= \E
[ \xi(\Po \cap (Q_\la)_+, i_1(\la^{1/d} x) + Q_{1})
\xi( \Po \cap (Q_\la)_+, i_1(\la^{1/d}x)+ z + Q_{1})  ]
\\
= \E [
\xi(\Po \cap (-i_1(\la^{1/d} x) + Q_\la)_+,  Q_{1})
\xi( \Po \cap (-i_1(\la^{1/d} x)  + Q_\la)_+,  z + Q_{1}) ] .
\eean
By a similar argument to \eq{0116c},
as $\la \to \infty$ we have
$$
\xi(\Po \cap (-i_1(\la^{1/d} x) + Q_\la)_+,  Q_{1})
\xi( \Po \cap (-i_1(\la^{1/d} x)  + Q_\la)_+,  z + Q_{1})
\toas
\xi(\Po ,  Q_{1})  \xi( \Po ,  z + Q_{1})
$$
and since $\xi$ is bounded \eq{bounded}, the expectations
converge. Hence, by \eq{0117c} and the similar lower bound,
 \bea
\E [\langle
 f_{\la,i_\la(x)} ,
 \xi_\la(\Po_\la \cap (Q_1)_+)
\rangle
\langle
 f_{\la,i_\la(x)+z} ,
\xi_\la(\Po_\la \cap (Q_1)_+) \rangle]
\nonumber \\
 \to f^2(x)\E[ \xi(\Po ,  Q_{1})  \xi( \Po ,  z +
Q_{1})].
 \lbl{0117d}
 \eea
 Also, $ \E [\langle
f_{\la,i_\la(x)}, \xi_\la(\Po_\la \cap (Q_1)_+) \rangle] $
converges to $
 f(x) \mu (\xi)
$
 by  \eq{0116a}, and a similar argument
yields
 \bean \E [\langle
 f_{\la,i_\la(x)+\la^{-1/d}z},
\xi_\la(\Po_\la \cap (Q_1)_+)
 \rangle ] \to f(x) \mu (\xi).
\eean Combining these with \eq{0117d}, we obtain that for $x \in
(0,1)^d$ with $f$ continuous at $x$, \bea \lim_{\la \to \infty}
\Cov \left[ \langle f_{\la,i_\la(x)}, \xi_\la(\Po_\la \cap (Q_1)_+
) \rangle, \langle f_{\la,i_\la(x)+z}, \xi_\la(\Po_\la \cap
(Q_1)_+ ) \rangle \right]
\nonumber \\
= f^2(x) \Cov \left[ \xi(\Po ,  Q_{1}) , \xi( \Po ,  z + Q_{1})
\right] 1_{Q_1}(x). \lbl{0123b} \eea Also, \eq{0123b}  holds for
$x \in \R^d \setminus [0,1]^d$ as well, since both sides are zero
for large $\la$. By \eq{expcov}, \eq{0123b}   and the dominated
convergence theorem, applied to  the last line of \eq{0116g}, we
obtain \bean
  \label{varlim2a} \lim_{\la \to \infty} \la^{-1} \Var [\langle f,
   \mu_\la^\xi \rangle]  = \sigma^2(\xi) \int_{[0,1]^d}
   f^2(x)dx.
\eean
In other words, we have demonstrated \eq{varlim2} in the
case where $f \equiv g$ and $f$ is nonnegative.
Extending \eq{varlim2} to the general case is then a routine
application of linearity.
%
%
$\qed$ \\

{\em Proof of \eq{rate} and the rest of Theorem \ref{basethm1}.}
Suppose $\sigma^2(\xi) > 0$ and take
$f \in \testfn $ with
  $\int_{Q_1} f^2(x)dx >0$.
We prove asymptotic normality
for $\langle f,\mu_{\la}^{\xi}\rangle$,
with a rate of convergence.
To do this we
  adapt the proof of Corollary 2.4 of
\cite{PY5}
(Corollary 2.1 in the electronically available version of \cite{PY5}),
 to the setting of functionals of TLF point sets in $\R^d$.
The proof of Corollary 2.4 of \cite{PY5} involves applying Stein's
method to a graph whose vertices are sub-cubes of the unit cube
with edge length proportional to $(\log \la) \la^{-1/d}$  and with
edges between sub-cubes whenever the distance between sub-cubes is
within twice the common cube edge length.
We make the
following trivial modifications to the proof of Corollary 2.4 of
\cite{PY5}.

Let $\la$ be fixed and large.  Subdivide
$Q_1^+$ into
 $V(\la):= 3^d
\la \r_{\la}^{-d}$ sub-cubes $C^{\la}_i$ of volume
$\la^{-1}\r_\la^d,$
where $\r_\la := \alpha \log \la$ for some suitably large
$\alpha$, as in section four of \cite{PY5}.
For all $1 \leq i \leq V(\la)$,
 put
$$
\xi^{\la}_{i}:=
 \int_{C_{i}^\la} f(x) \xi_\la(\Po_\la \cap Q_1, dx)
=
 \int_{\la^{1/d}C_{i}^\la} f(y) \xi(\Po \cap Q_\la, dy)
.
$$
Then
$$
\langle f,\mu_{\la}^{\xi} \rangle = \sum_{i=1}^{V(\la)}
\xi^{\la}_{i}.
$$
Note that $\xi_i^{\la}$ is the analog of $\sum_{j = 1}^{\infty}
|\xi_{ij} |$ of Lemma 4.3 of \cite{PY5} and furthermore, by the
 boundedness \eq{bounded} of $\xi$,
  for $q=3$ there exists $\Co:= \Co (q;f) < \infty$ such that
$||\xi_i^{\la} ||_q \leq \Co \r_{\la}^{d}$.

Consider for all $1 \leq i \leq V(\la)$ the events
$$
E_i := \bigcap_{j \in \Z^d: ( j + Q_1 ) \cap  \la^{1/d}
C^{\la}_i \neq \emptyset }
 \{R^{\xi}(j, \la) \leq \r_{\la} \},
$$
where $R^{\xi}(j,\lambda)$
is the radius of
stabilization of $\xi$ at $j \in  \Z^d$. Let
$$
E_\la := \bigcap_{i = 1}^{V(\la)} E_i,
$$
and note that $P[E_\la^c] \leq \la \tau(\r_{\la})$, where $\tau$
is as in Definition 2.1.

Next, define the analog of $T_\la'$ in \cite{PY5}
by
$$
{\mu'}_\la^{\xi}:= \sum_{i=1}^{V(\la)} \xi^{\la}_{i}{ \bf 1}_{E_i}
$$
and note that $\xi^{\la}_{i}{ \bf 1}_{E_i}$ and $\xi^{\la}_{j}{
\bf 1}_{E_j}$ are independent whenever $D_2(C_i^\la, C_j^\la) > 2
\la^{-1/d} \rho_\la$.
 For  $1 \leq i \leq V(\la)$,
 define
$$
S_i:= (\Var [{\mu'}_\la^{\xi}] )^{-1/2} \xi^{\la}_{i}{ \bf
1}_{E_i}
$$
and put $$ S := \sum_{i=1}^{V(\la)} (S_i - \E S_i).
$$

As in \cite{PY5} we define a dependency graph $G_\la := ({\cal
V}_\la, {\cal E}_\la)$ for $\{S_i\}_{i=1}^{V(\la)}$. The set
${\cal V}_\la$ consists of the sub-cubes
$C_1^\la,...,C_{V(\la)}^\la$ and edges $(C_i^\la, C_j^\la)$ belong
to $ {\cal E}_\la$ if $D_2(C_i^\la, C_j^\la) \leq 2 \la^{-1/d}
\r_\la.$  Next, in parallel with the proof of Corollary 2.4 of
\cite{PY5}, we notice that:
 \vskip.1cm (i)
$V(\la) := | {\cal V}_\la | = 3^d \la \rho_\la^{-d}$,

(ii) the maximal degree $D_\la$ of $G_\la$  satisfies $D_\la \leq
5^d$,

(iii) for all $1 \leq i \leq V(\la)$ we have
$ \| S_i\|_3 \leq \Co
(\Var [{\mu'}_\la^{\xi}] )^{-1/2} \rho_\la^{d}$,

(iv) $\Var [{\mu'}_\la^{\xi}] =
O( \rho_\la^{d} \la)$,

and

(v) $| \Var [ \langle f,\mu_{\la}^{\xi} \rangle]  - \Var [\langle
f, {\mu'}_\la^{\xi} \rangle] | \leq \Co \la^{-2}.$

\vskip.1cm

As in \cite{PY5}, we may use Stein's method to deduce a normal
approximation result for $S$ and then applying the estimates (iv)
and (v) and following \cite{PY5} verbatim we can turn this into a
normal approximation result for $\langle f,\mu_{\la}^{\xi}
\rangle$, i.e., in this way we obtain the
desired rate (\ref{rate}) when $f \equiv 1$. 

 The normal approximation result for
$\langle f,\mu_{\la}^{\xi} \rangle$, together with \eq{varlim2},
 implies that  $\la^{-1/2}
\langle f, \bar\mu_{\la}^{\xi} \rangle$ converges in
distribution to a mean zero
normal random variable with variance $\sigma^2(\xi) \int_{[0,1]^d}
f^2(x) dx.$
Given this,
the convergence of the
finite dimensional distributions in Theorem \ref{basethm1} is a
 standard application of the Cram\'er-Wold device.  This
completes the proof of Theorem  \ref{basethm1}.
\qed

\section{Stabilization of infinite input packing functionals}\label{STABII}
\label{secstab}
\allco
In this section, we show that the random packing measures $\nu_\la$
and $\nu'_\la$ described in Section \ref{INTRO} can each be expressed in
terms of a suitably defined measure-valued functional
$\xi$ of TLF point sets in $\R^{d} \times \R_+$,
 of the general type considered
 in Section \ref{TEAUX},
 applied
to a Poisson point process in space-time.
Then we show that in both cases the appropriate
choice  of $\xi$ satisfies the exponential stabilization condition
described in Definition \ref{stab}, so that Theorem \ref{basethm1}
is applicable to this choice of $\xi$. We defer to  the next
section the proof that in both cases the appropriate
 choice of $\xi$ satisfies
$\sigma^2(\xi) >0$.

Let us say that two points $(x,t)$ and $(y,u)$ in $\R^d \times \R_+$ are
{\em adjacent} if $(x+S) \cap (y+S) \neq \emptyset$.
  Given TLF 
$\X \subset \R^d \times \R_+$, let us first list the
  points of $\X$ in  order of  increasing time-marks
  using the  lexicographic ordering
  on $\R^d$
  as  a tie-breaker in the case of any pairs of points of $\X$
  with equal time-marks.
Then consider
the points of $\X$ in the order of the list; let the
first point in the list be accepted, and
 let each subsequent point be accepted if it
is  not
  adjacent to any previously  accepted point of $\X$;
otherwise let it be rejected.
  We call this the
  {\em usual rule} for packing points of $\X$, since
it corresponds to the packing rule of Section \ref{INTRO}
with the input ordering determined by time-marks.
Let ${\cal A}(\X)$ denote the subset of $\X$ consisting
of all accepted points when the points of $\X$ are packed
according to the usual rule.

We consider two specific measure-valued functionals $\xi^*$ and $\xi'$
on TLF point sets in $\R^d \times \R_+$,
of the general type
considered in Section \ref{TEAUX},
which are
 defined as follows.
 For any TLF
point set $\X \subset \R^d \times \R_+$
and bounded Borel $A \subset \R^d,$
recall that $A_+ := A \times \R_+$.
Let $\xi^*(\X,A)$ be the number of points of ${\cal A}(\X)$
which lie in  $A_+$, and with $|\cdot|$ denoting Lebesgue measure, let
$$
 \xi'(\X,A):= |S|^{-1}
\left| A \cap \left( \bigcup_{(x,t) \in {\cal A}(\X)}  (x+S)    \right)
\right|.
$$
Then $\xi^*$ and $\xi'$ are clearly translation
invariant, and  are bounded (i.e., satisfy \eq{bounded}),
since only a
bounded number of solids can be packed in any fixed bounded cube.

  Recall that $\P_\la $ denotes a homogeneous Poisson point process of
   intensity $\la$ on $\R^d \times \R_+$, and $\Po = \Po_1$. Assume
$\Po_\la$ is obtained from $\Po$ by $\Po_\la := \la^{-1/d} \Po$.
For all $\la >0 $, recall the definition of $\xi_\la$ in Section
\ref{TEAUX}, and  define the random measures
$$
{\mu}^{\xi^*}_{\la}(\  \cdot \ ) := \xi^*_{\la}({\cal P}_{\la}
\cap (Q_1)_+,\cdot) \ \  \text{and} \ \ {\mu}^{\xi'}_{\la} (\
\cdot \ ):= \xi'_{\la}({\cal P}_{\la} \cap (Q_1)_+,\cdot).
$$
Let $N_{\la}^{\xi^*}$ denote the total
mass of $\mu_\la^{\xi^*}$, i.e.
 \bean
N_{\la}^{\xi^*}:= \mu_{\la}^{\xi^*}({\cal P}_{\la} \cap [0,1]^d_+,[0,1]^d).
\eean
Then
$\mu^{\xi^*}_\la$
and
$\mu^{\xi'}_\la$
are the {\em  random packing point
 measure} and the
 {\em  random packing volume
 measure}, respectively, corresponding to
 the random sequential
adsorption process obtained by taking the spatial locations
of the points of  $\Po \cap Q_\la$, in order of increasing time-mark,
as the input sequence.  Since these spatial locations are independent
and uniformly distributed on $Q_\la$,
we have the distributional  equalities
\bea
 \label{totmass}
\mu^{\xi^*}_\la \eqd \nu_\la, \ \ \mu^{\xi'}_\la \eqd \nu'_\la,  \
\ \text{and} \ \
N^{\xi^*}_{\lad} \eqd N_{\lad},
\eea
where the measures $\nu_\la$ and $\nu'_\la$
are given  in (\ref{defnu})
and the jamming number $N_{\lad}$ is also given in
Section \ref{INTRO}.

We show in Lemmas \ref{Stablem1} and \ref{Stablem} below that both
$\xi^*$ and $\xi'$  are exponentially stabilizing, and therefore
we can apply  Theorem \ref{basethm1} to either of
 these choices of $\xi$.
To proceed with the proof of exponential stabilization,
consider a partition of  $\R^d$ into translates of
the unit cube $C:=Q_1 = [0,1)^d$.
%
%
  It is convenient to index these translates as
  $C_{i},\; i:= (i_1,\ldots,i_d) \in
  {\Z}^d,$ with $C_{i} := (i_1,\ldots,i_d) + C.$ We shall write
  $C_{i}^+ := \bigcup_{j \in {\Bbb Z}^d,\; ||i-j||_{\infty} \leq 1} C_j,$ that
  is to say $C_i^+$ is the union of $C_i$ and its neighboring cubes.
  We also consider the {\em moat}  $\Delta C_i := C_i^+ \setminus C_i.$


We need further terminology. Given TLF $\X \subset \R^d \times
\R_+$, and given $A \subset \R^d$, we say that $\X$ {\em fully
packs} the region $A$  if every point in $A_+$ is adjacent to at
least one  point of ${\cal A}(\X )$. For $t
>0$, we say $\X$ {\em fully packs $A$ by time t} if $\X \cap (\R^d
\times [0,t])$ fully packs $A$.
%
Given $B\subseteq \R^d$,
 we say that a  finite
 point configuration
  ${\cal X} \subset (B \cap C^+_{i})_+$ is {\em maximal} or
  {\it strongly saturates} the cube $C_{i}$ in $B$
  if for each TLF external configuration ${\cal Y}
\subset
  (B \setminus C^+_{i})_+,$
    ${\cal X} \cup {\cal Y}$ fully packs the region $B \cap C_{i}$
(the existence of  maximal configurations is guaranteed by Lemmas
\ref{lem3} and \ref{Probab}  below).

We shall be interested in strong saturation of $C_i$ in $B$
when $B= \R^d$ or when $B = Q_\la$. The reason for
our interest is this: If we
knew that there was a constant $\tau < \infty$ such that $\P \cap
(C_{\0}^+ \times [0,\tau])$ strongly
 saturated $C_{\0 } $ in $\R^d$ a.s., then
points in $\P$ with time marks exceeding $\tau$ would have no
bearing on the packing status of points in $\P \cap (C_{\0})_+$. Thus,
to check stabilization of $\xi$ at $\0$ it would be enough to
replace $\P$ by the
Poisson point process
$\P \cap (\R^d \times [0,\tau])$,
 and follow the stabilization arguments for packing with finite Poisson input (section four of
\cite{PY2}).  While clearly no such constant  $\tau $ exists,
we shall show in Lemma \ref{Probab} that  a finite random $\tau $ exists.



We say that $\X$ {\em locally strongly saturates} $C_i$
  if for each
  $ \eta \subseteq \X \cap (\Delta C_i)_+, $
  the point set $(\X \cap (C_i)_+) \cup \eta $
  fully packs  $C_i$.
 The following lemma shows that local strong saturation implies
strong saturation.


\begin{lemm}
\label{lem3}
 Suppose $\X \subset (C_i^+)_+$ is TLF and
  locally strongly saturates $C_i$.
Then
 for any $B \subseteq \R^d$ with
$C_i \subseteq B$,  $\X \cap  B$ strongly saturates $C_i$ in $B$.
 \end{lemm}
 {\em Proof.}
Let  $\Y \subset (B \setminus C_i^+)_+$ be TLF.
 Let  $\eta:= \A((\X \cap B_+) \cup \Y) \cap (\Delta C_i)_+$.
We claim  that
 \bea
\A((\X \cap B_+) \cup \Y) \cap (C_i^+)_+ =
 \A((\X\cap (C_i)_+)  \cup \eta) .
\lbl{0126}
\eea
 Indeed,  considering each point of $(\X \cap (C_i)_+) \cup \eta$
in the usual temporal order, we see that the decision on
whether to accept is the same for these points
 whether we are applying the usual packing
rule to $(\X \cap B_+) \cup \Y$ or to $(\X \cap (C_i)_+) \cup \eta$.

Since we assume $\X$
 locally strongly saturates $C_i$,
$(\X\cap (C_i)_+) \cup \eta$ fully packs $C_i$,
and so by \eq{0126},
$(\X \cap B_+) \cup \Y$
 fully packs $C_i$.
$\qed$ \\

We will use one more auxiliary lemma.

  \begin{lemm}\label{lem1}
  With probability 1, $\P$ has the property that
  for any
  $\eta \subseteq \P \cap (\Delta C_\0)_+$,  there exists
  $T < \infty$ such that the
  point set $(\P \cap (C_\0)_+) \cup \eta$
  fully  packs $C_\0$ by time $T$.
\end{lemm}
{\em Proof.}
Suppose that for
 each rational hypercube $Q$ contained in $C_\0$,
$\Po \cap Q_+ \neq \emptyset$;
this event has probability 1.

Take $\eta \subset \Po \cap (\Delta C_\0)_+$.
Let ${\cal A}:= \A ((\Po \cap (C_\0)_+) \cup \eta)$.
Clearly ${\cal A}$ is finite.
Let $ V$ be the set of $x \in  C_\0$ such that $(x,0)$
does not lie adjacent to any point of ${\cal A}$. Then
$V$  is open in $C_\0$
 (because we assume $S$ is closed) and if it is non-empty,
it contains a rational cube contained in $C_\0$ so that
$V_+$  contains
a point of $\Po \cap (C_\0)+$. But then this point should
have been accepted so there is a contradiction. Hence $V$
is empty and since ${\cal A}$ is finite this shows that $C_\0$
is fully packed within a finite time.
 $\qed$ \\

  For $i \in {\Z}^d$, let   $T_{i}:= T_{i}(\P)$ denote the
time  till local strong
  saturation,
  defined to be the smallest $t \in [0, \infty]$
  such that
$C_i $ is locally strongly saturated by
the point set  $(\P \cap (C^+_{i})_+) \cap (\R^d \times [0,t]) $
  (and set $T_i =\infty$ if no such $t$
exists).
  Clearly, $T_{i},\;
  i \in {\Z}^d,$ are identically distributed random variables depending only
  on $\P \cap (C^+_{i})_+$. In particular,
 $(T_i, \ i \in \Z^d)$ forms a 2-dependent random field,
 meaning that
 $T_i$ is independent of $(T_j, \|j-i\|_\infty >2)$
for each $i \in \Z^d$. We can now prove the key result that $T_\0$
is almost surely finite.

 \begin{lemm}\label{Probab} It is the case that $ P[T_{\0 } = \infty ] = 0.$
 \end{lemm}
{\em Proof.}
Suppose that $T_\0 = \infty$.
Then
 for each positive integer
$\tau$ there exists
$\eta_\tau \subseteq \P \cap (\Delta C_\0)_+$ such that
$(\P \cap (C_\0)_+) \cup \eta_\tau$
does not  fully pack $C_\0$
 by time $\tau$.

Assume
$ \P \cap (\Delta C_\0)_+$ is locally finite (this happens almost surely).
Then
$ \P \cap (\Delta C_\0 \times  [0,1])$ is finite so that
we can take a subsequence $\tau' \to \infty$
of $\tau$ along which
$\eta_{\tau'} \cap (\Delta C_\0 \times [0,1])$ is the same for
all $\tau'$. Then we can
take a further subsequence $\tau''$ of $\tau'$ along which
$\eta_{\tau''} \cap (\Delta C_\0 \times [0,2])$ is the same for
all $\tau''$. Repeating this procedure and using Cantor's diagonal
argument, we can find a subsequence $\tau_n$ tending to infinity,
and a limit set $\eta
\subset (\Delta C_\0 \times \R_+)$,
  such that for
all $k$, it is the case that
\bea
\eta_{\tau_n} \cap (\Delta C_\0 \times [0,k]) = \eta \cap
(\Delta C_\0 \times [0,k])
\label{1222}
\eea
 for all but finitely many $n$.


Let $k >0$, and choose $n$ to be large enough so that
$\tau_n \geq k$ and such that
(\ref{1222}) holds.
Then
the point set $(\P \cap (C_\0)_+) \cup \eta_{\tau_n} $ does not
 yet fully pack $C_\0$ by time $\tau_n$,
and therefore
 $(\P \cap (C_\0)_+) \cup \eta$ does not yet fully pack
$C_\0$ by time $k$.

Since $(\P \cap (C_\0)_+) \cup \eta$ does not yet fully pack
$C_\0$ by time $k$ for any $k$,
 we are in the complement of the event described in
Lemma \ref{lem1}. Thus by that result,
the event $\{T_\0 = \infty\}$ is contained in
an event of probability zero, which completes the proof of Lemma
\ref{Probab}. $\qed$ \\

Using Lemma \ref{Probab},
 we can now prove that $\xi^*$ and $\xi'$, defined at the start of this section,
satisfy the first part of  exponential stabilization
(exponential decay of the tail of $R'$).

\begin{lemm}\label{Stablem1}
  There exists a positive constant
${\Co}_1$ such that
for either $\xi = \xi^*$ or $\xi = \xi'$,
 there is a stabilization radius $R'$ as described
in Definition \ref{stab},  satisfying
  $$ P \left[R' > L\right] \leq \Co_1 \exp(-  L/K_1),
~~~~ \forall L >0. $$
 \end{lemm}

{\em Proof.}
 Let $\delta_1 > 0$ be a number falling below the critical probability
 $p_c$ for site percolation
 on ${\Z}^d$ with neighborhood relation ${i} = (i_1,\ldots,i_d) \sim {j} =
 (j_1,\ldots,j_d)$ if and only if
$ \|i-j\|_\infty \leq 1,$ see Grimmett \cite{Gr}.

 We will apply a {\it domination by product measures result} of \cite{LSS},
 more precisely Theorem 0.0 in \cite{LSS}. This tells us that, for
a family of $\{0,1\}$-valued random
 variables indexed by lattice vertices, if we are able to show that for each given site the probability
 of seeing $1$ there conditioned on the configuration outside a fixed size neighborhood of the site
 exceeds certain large enough $p,$ then this random field dominates a product measure with positive
 density $q$ which can be made
 arbitrarily close to $1$ by appropriate choice of $p$.
By this result, with $\delta_1 $ as chosen above
we can find $\delta_2 >0$
such that
any 2-dependent random field $(Y_i, i \in \Z^d)$ with $Y_i$
taking values in
 $\{0,1\}$ and $P[Y_i =1] \geq 1- \delta_2$ for each $i$, this random
field dominates the product measure with density $1 - \delta_1$.

Using Lemma \ref{Probab}, take
 $T^* > 0$ such that $P[T_\0 > T^*] < \delta_2$. Then
by the conclusion of the preceding paragraph,
  $\P$ can be coupled on a common probability space with an i.i.d.
 $\{0,1\}$-valued random field $\pi_i,\; i \in {\Bbb Z}^d,$ so that, for all $i \in {\Bbb Z}^d,$
 \begin{itemize}
  \item $P[\pi_i = 1] \geq 1-\delta_1,$
  \item we have $T_i < T^*$ whenever $\pi_i = 1.$
 \end{itemize}

Let us say that the cube $C_i$ is $T^*$-{\em saturated} if
$T_i \leq T^*$.
By Lemma \ref{lem3}, if $C_i$ is $T^*$-saturated then
for any $B \subseteq \R^d$ with $C_i \subseteq B$,
$\Po \cap ([C^+_i \cap B] \times [0,T^*])$ strongly saturates
$C_i$ in $ B$.

We declare
  a point $(x,t) \in \P \cap (C_{i})_+$ to be
{\em causally relevant} if either
   \begin{itemize}
    \item $\pi_i = 0,$
    \item or $\pi_i = 1$ and $t \leq T^*$.
   \end{itemize}
   Otherwise the point $x \in \P \cap (C_{i})_+$ is declared {\em causally irrelevant}.

   We now argue as follows, directly adapting the oriented percolation based technique
   introduced in section four of \cite{PY2}. We convert the collection of
points $\P$ (in $\R^d \times \R_+$)
   into a directed graph by providing a directed connection from $(y,s)$ to $(x,t)$ whenever
   $|y-x| \leq 2d_S$ and $s < t$ and, moreover, both $(x,t)$ and $(y,s)$ are causally relevant.
   By the {\em causal cluster} $\Cl[(x,t);\P]$ of $(x,t) \in \P$
we understand
   the set of all causally relevant
 points $(y,s)$ of $\P$ such that there is a directed path from
$(y,s)$ to  $(x,t)$
(referred to as a {\em causal chain for $(x,t)$} in the sequel).
   Necessarily the points in the causal cluster for $(x,t)$ have time mark
   at most $t$.

For each $(x,t) \in \P$ we define the {\em causal cube cluster} of
$(x,t)$ in $\R^d$
by
$$
\bar{  \Cl}[(x,t);\P] := \bigcup [C_j^+ : \ (C_j)_+ \cap \Cl[(x,t);\P] \neq
\emptyset ]
$$
and for each $i \in  \Z^d$ we define its {\em causal cube cluster
} as the union of clusters given by \be \label{ccc}\bar{
\Cl}[i;\P] := \bigcup_{(x,t) \in \P \cap (C_i^+)_+} \bar{
\Cl}[(x,t);\P].
 \ee

The significance of causal cube clusters is as follows. First,
we  assert that the packing status
of   a given point $(x,t)$ is unaffected by
changes to $\Po$ outside $\bar{\Cl}[(x,t);\Po]$.
Indeed, viewing
the directed connections as potential direct interactions between
overlapping solids in the course of the sequential packing
process, we can repeat the corresponding argument from
 Lemma 4.1 in \cite{PY2},
 adding
the extra observation  that causally irrelevant points will not be
accepted regardless of the outside packing
configuration and hence do not have to be taken into account.
Similarly, the packing status of the totality of points falling
within distance $d_S$ of the
cube $C_i$ can only be affected by the status of points falling in
the causal cube cluster $\bar{  \Cl}[i;\P]$. Consequently, we see
that for either
 $\xi = \xi^*$
or $\xi = \xi'$,
we can define a radius of stabilization
by
\be \label{Causal}
  R' :=
   \diam (\bar{ \Cl}[\0;\P]).
\ee
We need to show that $R'$ is almost surely finite with an
exponentially decaying tail. Given $L>0$, let $E_1(L)$ be the
event that there is a `path of zeros' from some site $i \in
\{-1,0,1\}^d$  to the complement of $B_{L/2 - \sqrt{d}}$ in the
Bernoulli random field $(\pi_i,i \in \Z^d)$. More formally,
$E_1(L)$  is the event that there exists there is  a sequence
$i_0, i_1, i_2, \ldots, i_n$, such that (a) $i_0 \in
\{-1,0,1\}^d$, and (b) $i_n \notin B_{L/2- 2 \sqrt{d}} $, and (c)
 for $j=1,\ldots,n$,
   $i_j \in \Z^d$ and $\|i_j -i_{j-1}\|_\infty =1$
 and $\pi_{i_j} = 0$.

For $i \in \Z^d$, let $E_2(L,i)$ be the event
that
there exists
$(x,t)\in \Po \cap (C_i)_+$, such that $t \leq T^*$ and there
exists a causal chain for $(x,t)$ which starts at some point of
$\Po \setminus (B_{L- 2 \sqrt{d}})_+$.
Define the event
$$
E_2(L) :=
\bigcup \left[
E_2(L,i):
i \in \Z^d, C_i  \cap B_{L/2} \neq \emptyset \right].
$$
Then we assert that
 the event $\{R' > L\}$ is contained in $E_1(L) \cup E_2(L)$.
Indeed, if $E_2(L)$ does not occur, then  for any causal
chain for any $(x,t) \in \Po \cap (C_\0^+)_+ $ starting outside
 $(B_{L- 2 \sqrt{d}})_+$, all points
in the causal chain of $(x,t)$ lying  inside  $(B_{L/2})_+$
must have time-coordinate greater than $T^*$;  if also $E_1(L)$
does not occur, at least one  of these points must lie in a cube
which is $T^*$-saturated, and therefore be causally irrelevant,
so in fact there is no causal chain for any
 $(x,t) \in \Po \cap (C_\0^+)_+ $ starting outside
$(B_{L- 2 \sqrt{d}})_+$.
Hence, $\bar{\Cl}[\0,\P] \subseteq B_L$, so that $R' \leq L$.

By the choice of $\delta_1$ and   by the exponential
decay of the cluster size in the subcritical percolation regime
(see e.g. Sections 5.2 and 6.3 in Grimmett \cite{Gr}), we have
exponential decay of $P[E_1(L)]$. That is,
there is a constant $\Co_2$ such that
 $P[E_1(L)] \leq \Co_2 \exp(-L/\Co_2 )$ for all $L$.

Since $T^*$ is fixed, we can
use the methods of \cite{PY2}
 for finite (Poisson) input packing, in particular
the argument leading to Lemma 4.2 in \cite{PY2}, to
see that  there is a constant $\Co_3$ such that
$P[E_2(L,i)] \leq \Co_3 \exp(-L/\Co_3)$
for all $i \in \Z^d \cap B_{L/2}$.
Since the number of such $i$ is only $O(L^d)$,
we see that
 $P[E_2(L)]$ also decays exponentially in $L$, and hence
so does $P[E_1(L)]+ P[ E_2(L)]$.
Since
 the event $\{R' > L\}$ is contained in $E_1(L) \cup E_2(L)$,
the lemma is proved. $\qed$ \\

To finish checking that $\xi^*$ and $\xi'$ satisfy the
 conditions for Theorem \ref{basethm1}, we consider
strong saturation, not only of unit cubes but of
cubes of slightly less than unit size. Let
$Q_{\zeta}^+$ denote the cube $[-\zeta^{1/d},2 \zeta^{1/d})^d$, i.e. the cube
of side $3 \zeta^{1/d}$ concentric with
$Q_\zeta$.
Let us say that $Q_\zeta$
is {\em locally strongly saturated}
 by a finite point set $\X \subset (Q_{\zeta}^+)_+$
if for every $\eta \subseteq \X \cap (Q_{\zeta}^+ \setminus Q_{\zeta})_+$,
 the point set
$(\X \cap (Q_{\zeta})_+) \cup \eta$ fully packs $Q_{\zeta}$.

\begin{lemm}
\lbl{unifsslem}
Given $\delta >0$, there exist constants $\eps >0$ and $t_0< \infty$ such that
for all $\zeta \in [1-\eps,1]$,
\bea
\lbl{0203}
P[ \Po \cap (Q_{\zeta}^+ \times [0,t_0]) ~ {\mbox{\rm  locally strongly
saturates }}
Q_\zeta ] > 1-\delta.
\eea
\end{lemm}
{\em Proof.}
By Lemma \ref{Probab}, we can choose $t_0$ such that
$\Po \cap (Q_1^+ \times [0,t_0])$ locally
strongly saturates $Q_1$,
 with probability at least $1 - \delta/2$.
Having chosen $t_0$ in this way, we can then choose $\eps$,
with $2d_S < (1-\eps)^{1/d}$,
so that for any $\zeta  \in [1-\eps,1]$,
\bean P[ \Po \cap ( (Q_1 \setminus Q_{\zeta}) \times [0,t_0] )
\neq \emptyset] < \delta/2.
\eean
For $\zeta <1$ with $2 d_S < \zeta^{1/d}$,
if
$\Po \cap (Q_1^+ \times [0,t_0])$
strongly saturates $Q_1$,
and
$ \Po \cap ( (Q_1 \setminus Q_{\zeta}) \times [0,t_0] ) $ is empty,
then $\Po \cap (Q_\zeta^+ \times [0,t_0])$ strongly saturates
$Q_\zeta$.  Hence, the preceding probability estimates complete
the proof. $\qed$

\begin{lemm}\label{Stablem}
  There exists a positive constant ${\Co}_4$
such that
for either $\xi = \xi^*$ or $\xi = \xi'$,
 there is a family  of stabilization radii $R(i, \la) = R^\xi(i,\la)$,
defined for $ \la \geq 1$ and
$i \in \Z^d $
as described in Definition \ref{stab},
which satisfy
  \bea
\sup_{\la \geq 1, i \in \Z^d
}
P \left[R(i,\la) > L\right] \leq \Co_4 \exp(-L/ \Co_4 ).
\lbl{0125}
 \eea
 \end{lemm}
{\em Proof.}  First let us restrict attention to $\la$ with
$\la^{1/d} \in \N$.
Adapting notation from the preceding proof, for
$(x,t) \in \Po \cap (Q_\la)_+$ we let $\Cl[(x,t);\Po \cap (Q_\la)_+]$
denote
   the set of all causally relevant
 points $(y,s)$ of $\P \cap (Q_\la)_+$ such that there is a directed path from
$(y,s)$ to  $(x,t)$, with all points in the path lying inside $(Q_\la)_+$.
Then
define the causal cube cluster in $Q_\la$ for
$(x,t)$ by
$$
 \bar{\Cl}[(x,t);\Po \cap (Q_\la)_+] := \bigcup [
C_j^+ \cap Q_\la: (C_j)_+ \cap
\Cl[(x,t); \Po \cap (Q_\la)_+] \neq \emptyset]
$$ and
and  for $i \in \Z^d $ by
$$
\bar{\Cl}[i;\Po \cap(Q_\la)_+]
= \bigcup_{(x,t) \in \Po \cap (Q_\la \cap C_i^+)_+ }
\bar{\Cl}[(x,t); \Po \cap (Q_\la)_+ ].
$$
Define
\bea
R(i,\la)
: =  \diam (\bar{\Cl}[i;\Po \cap (Q_\la)_+]), ~~~ \la^{1/d} \in \N.
\lbl{Causal2}
\eea
Then for $i \in  \Z^d$, the packing statuses of points
 of $\Po \cap (C_i^+ \cap Q_\la)_+$
are unaffected by changes to $\Po \cap (Q_\la)_+$ in
the region $(Q_\la \setminus B_{R(i,\la)}(i))_+$,
by the same argument as in the preceding proof.
Here we are using the fact that $\la^{1/d} \in \Z$,
and that
 if $C_i \subset Q_\la$ is $T^*$-saturated
then  $C_i$ is strongly saturated  in $Q_\la$ by
$\Po \cap  (Q_\la \times [0,T^*])$ (Lemma \ref{lem3}).
Thus,  $R(i,\la)$
 serves as a radius of stabilization
in the sense of Definition \ref{stab} (for either $\xi^*$ or $\xi'$).
Moreover,
$\bar{\Cl}[i;\Po \cap(Q_\la)_+] \subseteq
\bar{\Cl}[i;\Po ],
$
and so with $\Co_1$ as in the  the preceding proof we
have $P[R(i,\la) >L] \leq \Co_1 \exp(-L/\Co_1)$,
uniformly over $i, \la$ with $\la^{1/d} \in \N$.

Now suppose $\la^{1/d} \notin \N$. In this case,
instead of dividing $Q_\la$ into cubes of side 1, some of which
would not fit exactly, we divide $Q_\la$ into cubes of side slightly
less than 1, which do fit exactly, and repeat the above argument.

More precisely, we modify the proof of Lemma \ref{Stablem1}.
With $\delta_2$ as in that proof,
we  use Lemma \ref{unifsslem} to choose
constants
$\eps >0$  and $T^* < \infty$
(with $\max( 2d_S, 1/2) < (1-\eps)^{1/d}$)
in such a way that for any $\zeta \in [1-\eps,1]$ we have
\bean
P[ \Po \cap (Q_{\zeta}^+ \times [0,T^*]) ~ {\mbox{\rm  locally strongly saturates }}
Q_\zeta ] > 1-\delta_2.
\eean
With $\eps$ thus fixed, for all large enough $\la$ we can
choose $\zeta = \zeta (\la)  \in [1-\eps,1]$ in such a way that
$\la^{1/d}/\zeta^{1/d}$ is an integer.
Partitioning $\R^d$ into cubes $C'_i$ of volume  $\zeta$,
we can then follow the argument  already given for the
case $\la^{1/d} \in \N$, using the fact that the each of the
unit cubes $i+Q_1$,
for which we need to check conditions in Theorem \ref{basethm1},
is contained in the union of at most $2^d$ cubes in the partition
$\{C'_j\}$.
$\qed$

\section{Jamming variability, variance asymptotics}
\label{secJV}

\allco

At the end of this section,
 we complete the proofs of Theorems \ref{CLT1} and \ref{CLT2}.
First,
we
need
 to show that
the limiting variance $\sigma^2(S,d)$ is non-zero for all $d $ and all $S$.
This is achieved by
 Proposition \ref{CLT3}
and
Lemma \ref{VARIANCE}
below. The first of these results establishes that any convex
 $S \subseteq \R^d$ with
nonempty interior satisfies jamming variability (as defined in
remark 6, Section \ref{INTRO}),
 and the second establishes that this
is sufficient to guarantee that $\sigma^2 (S,d) >0$.
Recall from (\ref{totmass}) that we can work just as well with
$N_\la^{\xi^*}$ as with $N_\la$.

\begin{prop} \label{CLT3}
The convex body $S $
 has jamming variability.
\end{prop}

 {\em Proof.} Given $S$,
for all $x \in \R^d$
define
$$
\| x \| := \sup \{ a \geq 0 : (x + aS) \cap aS  = \emptyset\}.
$$
It is straightforward to verify that $\|\cdot\|$ is a norm on
$\R^d$, using the convexity of $S$ to verify the triangle
inequality. For nonempty $A \subset \R^d$, and $x \in \R^d$,
 write $D(x,A)$ for $\inf \{ \|x-y\|: y \in
A\}$. By our earlier assumption that $2 d_S < 1$ we have
 $\| x \| < \| x \|_\infty$ for all $x \in \R^d$.

For $\LL \subset \R^d$, we shall say $\LL$ is  {\em packed}
if $\|x-y\| \geq 1$ for all $x \in \LL$, $y \in \LL$,
and that $\LL $ is {\em maximally packed} if it is packed and
\bea D(w,\LL) < 1, ~~~ \forall w \in \R^d. \label{1220} \eea We
shall say $\LL$ is a {\em periodic set} if
for all $x \in \LL$ and $z \in \Z^d$ we have
$x + z \in \LL$.

Let $\LL$ be a maximally packed periodic subset of $\R^d$ (it is not hard
to see that such an $\LL$ exists).
Then the function $x \mapsto D(x,\LL)$ is a continuous function on
$\R^d$ that is periodic (i.e., $D(x,\LL) = D(x+z,\LL)$ for all $x
\in \R^d,z \in \Z^d$).  Hence the range of this function is the
continuous image of the compact torus $\R^d/\Z^d$, and so is
compact. Hence by (\ref{1220})  we have \bean \beta: =  \sup \{
D(w,\LL) : w \in \R^d \} < 1. \eean Then for $x \in \R^d$ and
$\alpha >0$, by scaling
\bea \lbl{1220a}
 D(x, \alpha \LL)
= \alpha
 D( \alpha^{-1}x,\LL) \leq \alpha \beta.
\eea
Choose $\delta >0$ such that $\beta(1+ 6 \delta) < 1 - 2
\delta$. For $i=1,2$, let $ \LL_i := (1+ 3i \delta) \LL. $
By (\ref{1220a}) and the choice of $\delta$ we have for all $x \in
\R^d$ and $i=1,2$  that
 \bea
 D(x,\LL_i) < 1- 2\delta.
\label{1220d} \eea

 Let $c_1$ denote the number of points
of $\LL$ in $[0,1)^d$. Denote by $\Bx(L)$ the hypercube
$[-L/2,L/2]^d$. For $i=1,2$, let $n_i(L)$ denote the number of
points of $\LL_i$ in $\Bx(L-4)$.
Then as $L \to \infty$, for $i=1,2$ we have \bea \label{1220b}
n_i(L) \sim c_1(1+3 \delta i)^{-d} L^d. \eea Let $n_3(L)$ denote
the maximum integer $m$ such that there exists a packed subset of
$\Bx(L) \setminus \Bx(L-6)$ with $m$ elements. Then there is a
finite constant $c_2$ such that for all $L \geq 6$ we have \bea
\label{1220c} n_3 (L) \leq c_2 L^{d-1}. \eea By (\ref{1220b}) and
(\ref{1220c}), we can choose $L_0$  such that for $L \geq L_0$ we
have 
\bea
  n_3(L) <  n_1(L) - n_2(L) . \label{1220e} 
\eea 
For $x \in
\R^d$ and $r>0$, set $\tB_r(x): = \{y \in \R^d: \|y-x\| \leq r\}$
(a ball of radius $r$ using the norm $\|\cdot\|$). For bounded $A
\subset \R^d$, let $T(A)$ denote the time of the first Poisson
arrival in $A$, i.e set
 $$
T(A) := \inf\{ t: {\cal P} \cap (A \times \{t\}) \neq \emptyset\},
$$
with the convention that the infimum of the empty set is $\infty$.
 Fix $L \geq L_0$, and for $i=1,2$  define the event
$E_i$ by
\bean
 E_i : = \left\{ \max  \{ T(\tB_\delta(x) ) : x \in
\LL_i \cap \Bx(L-4) \} < T \left( \Bx(L) \setminus \cup_{x \in
\LL_i \cap \Bx(L-4) } \tB_\delta(x) \right) \right\}.
\eean
 Let
 $i=1$ or $i=2$. If $y,y'$ are distinct
points of $\LL_i$ then $\|y-y'\| \geq 1+3\delta$. Hence, if also
$w \in \tB_\delta(y)$ and $w' \in \tB_\delta(y')$, then $\|w -w'\|
\geq 1+ \delta$ by the triangle inequality. Moreover, for $x \in
\R^d$, by (\ref{1220d}) and the triangle inequality  we can find
$y = y(x) \in \LL_i$ such that $\|x-w\|\leq 1-\delta$ for  all $w
\in \tB_\delta(y)$. Hence,
if $E_i $ occurs then
the set of accepted points (i.e., centroids of accepted shapes) of
the infinite input packing process on $\Bx(L)$ induced by $\Po$
  with arbitrary
external pre-packed configuration $\eta$ in $\R^d \setminus
\Bx(L)$, includes one point from each  $\tB_\delta(x), x \in \LL_i
\cap \Bx(L-4)$, and also contains no other points from $\Bx(L-6)$.

Thus for any pre-packed configuration $\eta$ in $\R^d \setminus
\Bx(L)$,
 if $E_1$ occurs the number
of accepted points in $\Bx(L)$ is at least $n_1(L)$, and if $E_2$
occurs the number of accepted points is at most $n_2(L) + n_3(L)$.
Also, the probabilities $P[E_1]$ and $P[E_2]$ are strictly
positive and do not depend on $\eta$. By (\ref{1220e}), it follows
that there is a constant $\epsilon>0$, independent of $\eta$, such
that ${\rm Var}[N_{L^d}^{\xi^*}(\Bx(L))|\eta] \geq \epsilon$. Thus
we have established the required jamming variability. \qed



 \begin{lemm}\label{VARIANCE}
It is the case that
  $ \liminf_{\la \to \infty} \la^{-1} \Var [N_{\lad}^{\xi^*}] >0.
 $
 \end{lemm}

 {\em Proof.}
By Proposition \ref{CLT3}, there exists $L>0$ such that
$\inf_{\eta} \Var N[[0,L]^d | \eta] > 0$, where the infimum is
over all admissible $\eta \subset \R^d \setminus [0,L]^d$.
We consider $\la$ with $\la^{1/d}/(L+4) \in \N$.
  We subdivide the cube $Q_\la$ into $n(\la) := \la/(L+4)^d$
equal-sized sub-cubes
  $\tilde{C}_{1,\la},\tilde{C}_{2,\la},\ldots,\tilde{C}_{n(\la),\la}$ arising
 as translates of $\Bx(L+4)$ centered at $x_{1,\la},\ldots,x_{n(\la),\la}$
respectively. For $1 \leq i\leq n(\la)$, let $\tilde{C}^-_{i,\la}$
be the translate of $\Bx(L)$ centered at $x_{i,\la}$, and let
$M_{i,\la}$ be the translate of
 $\Bx(L+2) \setminus
\Bx(L)$ centered at $x_{i,\la}$ (a `moat' around
$\tilde{C}^-_{i,\la}$).

Using terminology from Section \ref{secstab}, let $F_{i,\la}$ be
the event that the
 point set
$\Po \cap (M_{i,\la})_+$ fully packs
$M_{i,\la}$ by time 1, and let $G_{i,\la}$
be the event that $\Po \cap ((\tilde{C}_{i,\la}
\setminus M_{i,\la}) \times [0,1]) $ is empty.
Let $E_{i,\la} := F_{i,\la} \cap G_{i,\la}$. Then $ p := P[E_{i,\la}]$
 satisfies $p >0$, and does not depend on $i$ or $\la$.

   Observing that the events $E_{i,\la}, 1 \leq i \leq n(\la),$ are
 independent (the cubes  $\tilde{C}_i$
  are disjoint), denote the  (random) set of indices for which $E_{i,\la}$
   occurs by $I(\la) := \{i_1,...,i_{K(\la)} \}$.
 Then $\E [K(\la)] = p n(\la)$.
   Conditional on the event $E_{i,\la},$ the packing process inside $\tilde{C}_{i,\la}^-$
 has a particularly simple form - before time $1$ there are no
   points in $\tilde{C}_{i,\la}^-$, and after that time the newly arriving solids centered
   in $\tilde{C}_{i,\la}^-$ undergo the packing process according to the usual rules
  with the additional restriction that a solid overlapping  another one
   packed in $M_{i,\la}$ before time $1$ is rejected.
   Note that for $i \in I(\la)$, no new solids are accepted in $M_{i,\la}$
   after time $1$ and, moreover, the acceptance times
   of solids accepted in $M_{i,\la}$ before time $1$ have no
   influence on the behavior of the packing process in $\tilde{C}_{i,\la}^-$ after
time $1$;  only their spatial locations  matter.
   For a configuration $\eta$ of accepted points (only spatial locations
taken into account) in $M_{i,\la}$, the
   process described above will be referred to as packing in $\tilde{C}_{i,\la}^-$
   in the presence of the {\it pre-packed configuration} $\eta.$

Let ${\cal M}_\lambda$ be the sigma-algebra generated by the
points of $\Po \cap (Q_\lambda \times [0,1])$, i.e. the Poisson
arrivals up to time 1. Event $E_{i,\la}$ is ${\cal M}_\lambda$-measurable,
for each $i$.

By the conditional variance formula we have
  $$
   \Var\left[N^{\xi^*}_{\lad} \right] = {\Bbb E} \left[\Var\left( N^{\xi^*}_{\lad} | {\cal M}_{\la} \right)\right]
   + \Var \left[{\Bbb E}\left( N^{\xi^*}_{\lad} | {\cal M}_{\la} \right)\right]
   $$
   $$
   \geq
{\Bbb E} \left[ \left. \Var\left( N^{\xi^*}_{\lad} \right| {\cal M}_{\la}
\right)\right]
 $$
$$
= \E \Var \left[ \left. \sum_{i \in I(\la)}
N^{\xi^*}_{\lad}[\tilde{C}^-_{i,\la}] + \left(  N_{\lad}^{\xi^*} -
\sum_{i \in I(\la)} N^{\xi^*}_{\lad}[\tilde{C}^-_{i,\la}] \right) \
\right| \  {\cal M}_{\la} \right],
$$
where we set $N^{\xi^*}_{\lad} [\tilde{C}^-_{k,\la}] := \xi^*(\Po
\cap Q_\la, \tilde{C}^-_{k, \la})$,  the number of solids packed in
$\tilde{C}^-_{k,\la}.$ Conditionally on ${\cal M}_{\la}$, the
packing processes after time 1 over different sub-cubes
$\tilde{C}^-_{i,\la}, i \in I(\la)$, are independent of each other
and of the packing process after time 1 in $Q_\la \setminus
\cup_{i \in I(\la)} \tilde{C}^-_{i,\la}$. Hence,
  \bean
   \Var\left[ N^{\xi^*}_{\lad} \right] \geq
    \E \sum_{i \in I(\la)} \left[ \Var [  N^{\xi^*}_{\lad}[\tilde{C}^-_{i,\la}] \ | \
    {\cal M}_{\la} ] \right] \geq \E [K] \inf_{\eta} \Var N[[0,L]^d|\eta],
\eean
  where the infimum is taken over all admissible configurations $\eta$
 outside $[0,L]^d$,
  and where $N[[0,L]^d|\eta]$ stands  for number of solids packed in
$[0,L]^d$ in the
  presence of the pre-packed configuration $\eta.$
By  Proposition \ref{CLT3}, this infimum is strictly positive, and
Lemma \ref{VARIANCE} follows.
\qed \\

 {\em Proof of Theorems \ref{CLT1} and \ref{CLT2}.}
Let $\xi$ be  $\xi^*$  as defined in Section \ref{secstab}.
Then
 Lemmas \ref{Stablem1} and \ref{Stablem} show that
$\xi=\xi^*$ satisfies the exponential stabilization
conditions in  Theorem \ref{basethm1}, so it satisfies
the conclusions \eq{explim2}, \eq{varlim2} and \eq{explim}
of that result.   The conclusion \eq{explim} gives us
\eq{meanrate} of Theorem \ref{CLT1}. Also, by putting
$f \equiv g \equiv 1$ on $Q_1^+$ and using \eq{varlim2}, we obtain the
variance convergence $\la^{-1}\Var N_\la \to \sigma^2$ asserted
in Theorem \ref{CLT1}.
By Lemma \ref{VARIANCE}, we may therefore deduce that $\sigma^2>0$.
Hence we may apply the last part of Theorem \ref{basethm1}
to obtain the rest of the conclusions in Theorem \ref{CLT2}
as they pertain to $\nu_\la$; also the conclusion
\eq{rate} of Theorem \ref{basethm1}
gives us \eq{jamrate}.

To get the same results for $\nu'$, we argue similarly with $\xi =
\xi'$. We need to check that the limiting means and variances are
the same, i.e. $\mu(\xi') = \mu(\xi^*)$ and $\sigma^2(\xi') =
\sigma^2(\xi^*)$. To see this, note that if  $f  \equiv 1$ on
$Q_1^+$, then $ \langle f, \mu_\la^{\xi'} \rangle = \langle f,
\mu_\la^{\xi^*} \rangle
$
so  application of \eq{explim2}  to this choice of $f$
yields
 $$
\mu(\xi') =
\lim_{\la \to \infty} \la^{-1} \E [\langle f,
  \mu_\la^{\xi'} \rangle]
=
\lim_{\la \to \infty} \la^{-1} \E [\langle f,
  \mu_\la^{\xi^*} \rangle]
=
\mu(\xi^*)
$$
and a similar argument using \eq{varlim2} shows that
 $\sigma^2(\xi') = \sigma^2(\xi^*)$.
  $\qed$

\vskip.5cm

Tomasz Schreiber, Faculty of Mathematics and Computer Science,
Nicholas Copernicus University, Toru\'n, Poland: \ {\texttt
tomeks@mat.uni.torun.pl }

 \vskip.5cm

Mathew D. Penrose, Department of Mathematical Sciences, University
of Bath, Bath BA2 7AY, United Kingdom: {\texttt
m.d.penrose@bath.ac.uk}

\vskip.5cm

J. E. Yukich, Department of Mathematics, Lehigh University,
Bethlehem PA 18015, USA:

{\texttt joseph.yukich@lehigh.edu}

\end{document}